\documentclass[12pt]{amsart}
\usepackage{amscd,amsmath,amssymb,amsfonts}
\theoremstyle{plain}
\newtheorem{thm}{Theorem}

\newtheorem{cor}[thm]{Corollary}
\newtheorem{prop}[thm]{Proposition}

\theoremstyle{definition}

\newtheorem{rmks}[thm]{Remarks}

\numberwithin{thm}{section}
\numberwithin{equation}{section}

\newcommand{\eq}[2]{\begin{equation}\label{#1}#2 \end{equation}}
\newcommand{\ml}[2]{\begin{multline}\label{#1}#2 \end{multline}}
\newcommand{\ga}[2]{\begin{gather}\label{#1}#2 \end{gather}}

\newcommand{\surj}{\twoheadrightarrow}

\newcommand{\sO}{{\mathcal O}}

\newcommand{\C}{{\mathbb C}}

\newcommand{\F}{{\mathbb F}}

\renewcommand{\H}{{\mathbb H}}

\newcommand{\N}{{\mathbb N}}
\renewcommand{\P}{{\mathbb P}}
\newcommand{\Q}{{\mathbb Q}}

\newcommand{\Z}{{\mathbb Z}}

\begin{document}

\title[Fano varieties]{Congruences for the
number of rational points, Hodge type and motivic conjectures for Fano varieties}
\author{Spencer Bloch}
\address{University of Chicago, Mathematics, IL 60 636, Chicago, USA}
\email{bloch@math.uchicago.edu}
\author{H\'el\`ene Esnault}
\address{Mathematik,
Universit\"at Essen, FB6, Mathematik, 45117 Essen, Germany}
\email{esnault@uni-essen.de}
\date{Dec. 17, 2002}
\begin{abstract}
A Fano variety is a smooth, geometrically connected variety over a
field, for which the dualizing sheaf is anti-ample. For example
the projective space, more generally flag varieties are Fano
varieties, as well as hypersurfaces of degree $d\le n$ in $\P^n$.
We discuss the existence and number of rational points over a
finite field, the Hodge type over the complex numbers, and the
motivic conjectures which are controlling those invariants. We
present a geometric version of it.
\end{abstract}
%\subjclass{Primary Algebraic Geometry}
\maketitle
\begin{quote}

\end{quote}

\section{Congruence for the
number of rational points for a variety over a finite field}

Let $X$ be a smooth projective variety over a field $k$. If $k=\F_q$
is finite,
it will be rarely the case that $X$ has a rational point. Yet, if the variety
is very {\it negative} in the sense of differential geometry, then
$X$ will have some, or even many rational points. The simplest example is
the projective space $\P^n$. Since one description of the
$k$-rational points $|\P^n(k)|$
is the quotient of the punctured vector space
$k^{n+1}\setminus\{0\}$ by the diagonal action of the homotheties $k^\times$,
 one sees that
\ga{1}{|\P^n(\F_q)|=\frac{q^{n+1}-1}{q-1}=1+q +\ldots +q^n.}
One way to measure  $|X(\F_q)|$ for a smooth projective variety
$X\subset \P^n$ is to consider the congruence
\ga{2}{|X(\F_q)|\equiv |\P^n(\F_q)| \ \text{mod}  \ q^\kappa,}
or equivalently
\ga{3}{|U(\F_q)|\equiv 0 \ \text{mod} \ q^\kappa,}
where $U=\P^n\setminus X$,
for some natural number $\kappa$. Of course if $\kappa=0$, then \eqref{2}
says nothing, but if there is a $\kappa \ge 1$ which satisfies
\eqref{2}, then first of all,
$X$ has a rational point,
and secondly the larger $\kappa$, the more rational points.
One codes $|U(\F_{q^s})|$ for all finite extensions $\F_q\subset \F_{q^s}$
in the zeta function defined by its logarithmic derivative
\ga{4}{\frac{\zeta'(U,t)}{\zeta(U,t)}=\sum_{s\ge 1} |U(\F_{q^s})|t^{s-1}.}
Thus the existence of a $\kappa \in \N\setminus \{0\}$ as in \eqref{3} for
 all finite extensions of $\F_q$
is equivalent to $\zeta(U,t)$, as a power expansion
in $t$, be in $\Z[[q^\kappa t]]$.
On the other hand, the fundamental theorem of Dwork \cite{Dw} asserts that
$\zeta(U,t)$ is a rational function over the rational numbers
\ga{5}{\zeta(U,t)\in \Q(t).}
One concludes that writing
\ga{6}{\zeta(U,t)=\frac{\prod_{i=1}^a (1-\alpha_it)}{\prod_{j=1}^b
(1-\beta_j t)},}
the reciprocal roots $\alpha_i$  and poles $\beta_j$
of the $\zeta(U,t)$ are divisible
by $q^\kappa$ as  algebraic integers, i.e. in $\bar{\Z}\subset \bar{\Q}$.
On the other hand, the Grothendieck-Lefschetz fixed point
trace formula \cite{Gr} asserts
\ga{7}{\zeta(U,t)=\prod_{i=0}^{2\ \rm dim(U)} {\rm det}(1-F_it)^{(-1)^{i+1}}}
where $F_i$ is the arithmetic Frobenius acting on the compactly supported
$\ell$-adic cohomology  $H^i_c(\bar{U}, \Q_\ell)$, which is isomorphic
to the primitive cohomology $H^{i-1}_{{\rm prim}}(\bar{X},\Q_\ell)=
H^{i-1}(\bar{X}, \Q_\ell)/H^{i-1}(\bar{\P^n},\Q_\ell)$ of $X$
for $(i-1)\le 2 \ {\rm dim}(X)$, and $=H^i(\bar{\P^n}, \Q_\ell)$ for $i\ge 2 {\rm dim}(X)$.
By the Weil conjectures proven by Deligne \cite{DeWeI}, the eigenvalues
of $F_i$ in any complex
embedding $\Q_\ell \subset \C$ have absolute values $q^{\frac{(i-1)}{2}}$.
Thus in \eqref{7} there can't be any cancellation between odd and even $i$'s.
Our condition \eqref{3} translates then exactly into the condition that
the eigenvalues of $F_i$ be all divisible by $q^\kappa$ as algebraic integers.

On the other hand, the cohomology $H^m(\bar{X}, \Q_\ell)$ carries
a {\it coniveau filtration} $\ldots \subset N^a H^m(\bar{X},
\Q_\ell)\subset N^{a-1} H^m(\bar{X}, \Q_\ell)\subset \ldots $,
where $N^a$ is the subgroup of classes which die after restriction
outside of  a codimension $a$ subscheme. One checks by a
d\'evissage as in \cite{E2}, Lemma 2.1, that the eigenvalues of
the Frobenius acting on $ N^aH^m(\bar{X}, \Q_\ell)$ are divisible
by $q^a$ as algebraic integers. Thus if $N^\kappa H^m_{{\rm
prim}}(\bar{X}, \Q_\ell)= H^m_{{\rm prim}}(\bar{X}, \Q_\ell)$,
the congruence \eqref{3} holds. The  Tate conjecture predicts the
converse: if the eigenvalues of Frobenius acting on $H^m_{{\rm
prim}}(\bar{X}, \Q_\ell)$ are all divisible as algebraic integers
by $q^\kappa$, then the primitive cohomology should be supported
in codimension $\kappa$, that is in concrete terms, there should
exist a codimension $\kappa$ subscheme $Z\subset X$ such that the
restriction map \ga{8}{H^m_{{\rm prim}}(\bar{X}, \Q_\ell) \to
H^m(\overline{X\setminus Z}, \Q_\ell)/H^m(\bar{\P^n}, \Q_\ell)}
dies.

\section{Hodge type} \label{hodgetype}
The {\it Hodge type} of a projective variety $X\subset \P^n$ over
a field $k$ of characteristic 0 is defined to be the largest
natural number $\kappa$ such that the Hodge filtration $\ldots
\subset F^aH^m_{DR, {\rm prim}}(X)\subset F^{a-1} H^m_{DR, {\rm
prim}}(X)\subset \ldots$ fulfills \ga{9}{F^\kappa H^m_{DR, {\rm
prim}}(X)=H^m_{DR, {\rm prim}}(X)} for all $m$. (One would define
similarly the Hodge type of $H^m_{DR, {\rm prim}}(X)$ for a given
degree $m$). If $X$ is projective smooth, then Hodge type
$=\kappa$ means \ga{10}{H^q(X, \Omega^p_X)=\begin{cases}
0 \ \text{for} \ q\neq p \ <\kappa\\
k \ \text{for}  \ q= p< \kappa.\end{cases}\notag}
If
\ga{10b}{N^\kappa H^m_{DR, {\rm prim}}(X)= H^m_{DR, {\rm prim}}(X),}
with
the same definition of the coniveau filtration for de Rham cohomology
as for $\ell$-adic cohomology,  then
 one easily computes via the Gysin sequence
 that the Hodge type of the primitive
cohomology is $\ge \kappa$. The Hodge conjecture predicts the converse:
if the Hodge type is $\kappa$, then the primitive
cohomology is supported in codimension $\kappa$.

\section{Hodge cohomology and slopes}
In this section,  $X$ is still assumed to be a smooth projective
variety over a finite field $\F_q$, where $q=p^d$ for a prime
number $p$. We are interested in conditions which force the
reciprocal zeros $\alpha_i$ and poles $\beta_j$ of  the zeta
function \eqref{6} to be divisible by $q^\kappa$ as algebraic
integers.

Write $W(\F_q)$ for the ring of Witt vectors over $\F_q$. It is a
complete discrete valuation ring with residue field $\F_q$. For
example, $W(\F_p) = \Z_p$, the $p$-adic integers. Let $K$ be the
quotient field of $W(\F_q)$. Frobenius on $\F_q$ induces
automorphisms $\sigma$ of $W(\F_q)$ and $K$ satisfying
$\sigma^d=\text{identity}$. The crystalline cohomology \cite{Ber}
$H^*_{{\rm crys}}(X/W(\F_q))\otimes K$ is a finite dimensional
$K$-vector space with an endomorphism $f$ (Frobenius) satisfying
 \eq{s2}{f(wx) = \sigma(w)f(x);\quad w \in K,\ x \in H^*. } In
particular, $f^d$ is $K$-linear, and the basic theorem \cite{Ber}
is that
 \eq{s3}{\zeta(X/\F_q,T) = \det(1-f^dT|H^*)^{-1}, } where of
course the determinant of the right is taken in the graded sense,
i.e. characteristic polynomials coming from $H^{{\rm odd}}$ appear
in the numerator of $\zeta$.

In order to estimate divisibility for the eigenvalues, we can
calculate crystalline cohomology using the de Rham-Witt complex
\cite{I}. This is a complex of pro-sheaves for the Zariski
topology on $X$
 \eq{s4}{W_\bullet\Omega^* := \{W_\bullet \sO
\xrightarrow{d} W_\bullet \Omega^1 \xrightarrow{d} W_\bullet
\Omega^1 \to \cdots \xrightarrow{d} W_\bullet \Omega^{\dim X}\}. }
In other words, each $W_\bullet \Omega^i$ is a projective system
of Zariski sheaves on $X$
 \eq{s5}{\ldots\surj W_n\Omega^i\surj
W_{n-1}\Omega^i\surj \ldots \surj W_1\Omega^i := \Omega^i_X }
where $\Omega^i_X$ is the sheaf of K\"ahler $i$-forms on $X$. Each
$W_n\Omega^i$ has a finite filtration with graded pieces
coherent, so the cohomology groups $H^j(X,W_n\Omega^i)$ have
finite length. For $i=0$, $W_n\sO$ is the sheaf of Witt vectors of
length $n$ over the structure sheaf $\sO_X$.

It is not true in general that the $H^j(X,W_\bullet \Omega^i):=
\varprojlim_n H^j(X,W_n \Omega^i)$ are finitely generated
$W(\F_q)$-modules (even for $i=0$,  \cite{Serre}). However, the
groups $H^j(X,W_\bullet \Omega^i)\Big/(\text{torsion})$ are
finitely generated. In particular, the $H^j(X,W_\bullet
\Omega^i)\otimes K$ are finite $K$-vector spaces. The
differentials in \eqref{s4} come from differentials $W_n\Omega^i
\xrightarrow{d} W_n\Omega^{i+1}$, and we define
 \eq{s6}{\H^*(X, W_\bullet\Omega^*)
:= \varprojlim_n \H^*(X, W_n\Omega^*). }

The de Rham-Witt complex plays the role of a sort of de Rham
complex calculating crystalline cohomology. Namely, there is a
canonical, functorial isomorphism
 \eq{s7}{H^*_{{\rm crys}}(X/W(\F_q))
\cong \H^*(X, W_\bullet\Omega^*). } Crucial for our purposes is
that the frobenius $f$ has a nice description on the de Rham-Witt
cohomology. Namely, one has endomorphisms
 \eq{s8}{f^{(i)}, v^{(i)}
: W_\bullet \Omega^i \to W_\bullet \Omega^i } which satisfy
$f^{(i)}v^{(i)} = v^{(i)}f^{(i)} = p$. The endomorphism $v^{(i)}$
is topologically nilpotent in the sense of the inverse system
\eqref{s5}. One has $d\circ f^{(i)} = pf^{(i+1)}\circ d$, so in
particular the $p^if^{(i)}$ on $W_\bullet \Omega^i$ induce a map
of complexes on $W_\bullet\Omega^*$. The resulting map on
$\H^*(X, W_\bullet\Omega^*)$ coincides with the Frobenius $f$ on
crystalline cohomology under the isomorphism \eqref{s7}.

As a consequence of these facts, one deduces that the spectral
sequence
 \eq{9}{E_1^{ab} = H^b(X, W_\bullet\Omega^a)\otimes K
\Rightarrow H^{a+b}(X, W_\bullet\Omega^*)\otimes K } degenerates
at $E_1$ (\cite{Bcrys}), and that the eigenvalues of the
$K$-linear endomorphism $q^af^{(a)d}:H^*(X, W_\bullet\Omega^a)
\to H^*(X, W_\bullet\Omega^a)$ coincide with the eigenvalues
$\alpha_i$ and $\beta_j$  appearing in $\zeta(X/\F_q,T)$ which
are divisible by $q^a$ but not by $q^{a+1}$. (This is because
$q^av^{(a)d}f^{(a)d} = q^ap^{d} = q^{a+1}$, and $v^{(a)}$ is
topologically nilpotent.)

Of course, there is a lot of mathematics here, and it is not
possible to give the details in a survey such as this. Note,
however, that there are two deep global results, \eqref{s3} and
\eqref{s7}. The rest involves the definition and local structure
of $W_\bullet\Omega^*$.

As a corollary of the above, we deduce
\begin{cor}\label{cor1} Let $\kappa\ge 1$ be a given integer.
Then all the reciprocal zeroes and poles $\alpha_i$ and $\beta_j$
of $\zeta(X/\F_q,T)$ are divisible by $q^\kappa$ if and only if
$H^*(X,W_\bullet\Omega^a)\otimes K = (0)$ for $a<\kappa$.
\end{cor}

We would like a criterion in terms of the Hodge groups
$H^b(X,\Omega^a_X)$ which will insure the de Rham-Witt groups
vanish as in corollary \ref{cor1}. The following is deduced from
a purely local calculation using the structure of the sheaves
$W_n\Omega^i$. We use the notation ``Hodge type'' as in section
\ref{hodgetype}, even though the ground field is finite.

\begin{prop} \label{slopes} With notation as above, if $X$ has Hodge type $
\kappa\ge 1$, then $H^*(X, W_\bullet\Omega^a) = (0)$ for
$a<\kappa$. In particular, all the reciprocal zeroes and poles of
$\zeta(X/\F_q,T)$ are divisible by $q^\kappa$.
\end{prop}
\begin{proof} We will use some results about the structure of $W_\bullet\Omega^*$ from
\cite{Ill}. The first point is that $f^{(i)}, v^{(i)}, p$ are
injective on pro-objects. Topological nilpotence for $v^{(i)}$
means
 \eq{s10}{v^{(i)n}W\Omega^i \subset \ker
(W_\bullet \Omega^i \to W_n\Omega^i) } \cite{Ill}, (2.2.1).
Further, by op. cit. (2.5.2), there is an exact sequence
\eq{s11}{0 \to W_\bullet\Omega^{i-1}/f^{(i-1)}
\xrightarrow{dv^{(i-1)}} W_\bullet\Omega^i/v^{(i)} \to \Omega^i_X
\to 0. } (For $i=0$ this says $W/v^{(0)} \cong \sO_X$.) By
induction on $i$ we see that $H^*(X,W_\bullet\Omega^i/v^{(i)}) =
(0)$ for $i<k$, so \ml{12}{v^{(i)n}: H^*(X, W\Omega^i) \cong
H^*(X,
W\Omega^i) = \\
\varprojlim H^*(X, W_n\Omega^i) \subset \prod_n H^*(X,
W_n\Omega^i). } By
 \eqref{s10}, we deduce that $H^*(X,
W_\bullet\Omega^i) = (0)$ for $i<k$.
\end{proof}
\begin{rmks} \label{rmk:ax}
\begin{itemize}
\item[i)] Proposition \ref{slopes} is a special case of a more
general theorem asserting that the Newton polygon of the
$F$-crystal $H^m(X/W)/({\rm torsion})$ lies above the Hodge
polygon defined by slope $i$ with multiplicity ${\dim }\
H^{m-i}(X, \Omega^i_X)$ (see \cite{Ma}, \cite{Ogus}). We have seen
that only the local properties of $W_\bullet \Omega^*$ are
relevant for Proposition \ref{slopes}.
\item[ii)]
It is of course very easy to compute the Hodge cohomology
$H^j(X,\Omega^i)$ for smooth complete intersections $X\subset
\P^n$ defined by $r$ equations of degrees $d_1\ge d_2\ge \ldots
\ge d_r$. It is $\kappa=[\frac{n-d_2-\ldots -d_r}{d_1}].$ Here
$[z]$ is the integral part of the rational number $z$.  In other
words, Proposition \ref{slopes} is an easy proof of the theorem
of Ax and Katz (\cite{Ka}) asserting \eqref{2} or equivalently
\eqref{3} in this case.
\end{itemize}

\end{rmks}

\section{From $\F_q$ to $\C$ and vice-versa.}
Let us think now that our smooth  variety over a finite field is coming
via reduction modulo $p$ from a variety defined in characteritic 0,
over a ring of finite type over the integers. Then, via
the comparison between $\ell$-adic and de Rham cohomologies,
 the coniveau in which
(primitive) de Rham cohomology is carried is the same as the coniveau in which
(primitive) $\ell$-adic is carried. In conclusion, one sees that
 the coincidence
of the $\kappa$ stemming from the $\zeta$ function with the $\kappa$ from
the Hodge type is a test both for the Tate and the Hodge conjectures.

We have two tests at disposal. First {\it smooth complete
intersections}. Let $X\subset \P^n$ be a smooth complete
intersection defined by $r$ equations of degree $d_1\ge d_2\ldots
\ge d_r$. We define
 $\kappa =[\frac{n-d_2-\ldots -d_r}{d_1}]$
 as in Remark
 \ref{rmk:ax}.
 Then, as already mentioned,
 the theorem of Ax and Katz \cite{Ka} asserts \eqref{3}
while Deligne's theorem \cite{DeSGA} asserts \eqref{9}. Moreover,
this bound is sharp both on the $\zeta$  and on the Hodge sides.

The smooth complete
intersections just discussed  with $\kappa\ge 1$ are special Fano
varieties.  We consider now our second example: {\it Fano varieties}
which are abstractly defined.
In characteristic 0, Kodaira vanishing applied to the ample
invertible sheaf $\omega^{-1}$ yields $H^q(X, \sO_X)=0$ for all $q>0$.
Thus the Hodge theoritic $\kappa$ is at least 1.
Over a finite field, \cite{E2}, Corollary 1.3 asserts that $|X(\F_q)|
\equiv 1 \ \text{mod}\ q$. Thus the  $\kappa$ of the $\zeta$ function
is at least 1 as well. This is our second test.
However, we observe that the test is not complete.
It might well be that the Hodge type of $X$ is $\ge 2$. Yet
 the proof given in
\cite{E2} does not give a better
the congruence for the number of rational points
over a finite field, unless we know the Chow groups of $X$.
This is the subject of the next section.

\section{Motivic conjectures}
The Beilinson-Bloch conjectures (\cite{B}, \cite{Be1}, \cite{Be2})
predict that the Chow groups of  a  smooth projective variety
defined over the complex numbers should be controlled by its
Hodge theory. More precisely, it predicts that if the Hodge type
of the primitive cohomology of a smooth complex projective
variety $X$ of dimension $d$ is $\kappa$, that is \eqref{9} holds
true, then one has \ga{11}{ CH_i(X)\otimes \Q=H^{2(d-i)}(X_{{\rm
an}}, \Q) \ \text{for} \ 0\le i\le (\kappa-1).} Applying the
splitting of the diagonal as initiated in \cite{B}, Appendix to
Lecture 1, and then refined as in \cite{J}, \cite{P}, \cite{EL},
one sees that \eqref{11} is equivalent to saying that there is a
nontrivial natural number $N$, there are dimension $i$ and
codimension $i$ cycles $\alpha_i$ and $\beta^i \subset X$,
 a codimension $\kappa$ cycle $Z\subset X$, a codimension $d$
cycle $\Gamma\subset X\times Z$ with
\ga{12}{N\Delta \equiv \alpha_0\times X + \alpha_1\times \beta^1+\ldots
+\alpha_{\kappa -1}\times \beta^{\kappa -1} + \Gamma}
where $\equiv$ means the equivalence
in $CH^d(X\times X)$, and $\Delta$ is the diagonal.
It is easily seen that such a decomposition \eqref{12} of the
diagonal,  up to torsion,
 implies \eqref{10b}. Over the complex numbers,
 it implies  a fortiori \eqref{9} while over a finite field,
 it implies \eqref{3}, by showing
again that the eigenvalues of the Frobenius are
 divisible by $q^\kappa$ as algebraic numbers.

In other words, the Beilinson-Bloch conjectures on the Chow groups
imply here the Hodge conjecture. Thus we expect that our Fano
varieties with Hodge type $\kappa$ will fulfill \eqref{11} and
\eqref{3} as well.

Let us look at our two classes of examples. The smooth complete intersections
as in section 3 with $\kappa \ge 1$ fulfill $CH_0(X)=\Z$ by a theorem
of Roitman \cite{Ro}. However, we do not know in general whether
 $CH_i(X)\otimes \Q=\Q$ for $i<\kappa$. We have very  few simple examples
where the bound is achieved (\cite{V}, \cite{O}), and the general results
 we have yield  bad bounds
( \cite{ELV}). We observe nevertheless that Roitman's theorem
(loc. cit.) yields a whole class of examples for which the
Beilinson-Bloch conjecture is true and sharp. Let $Y\subset
\P^{\kappa -1} \times \P^n$ be a smooth complete intersection of
bidegree $(1,d)$. Then one easily computes that
\ga{12a}{ H^q(Y, \Omega^p)=0 \ \text{for} \ q \neq  p, p \le (\kappa -1)\\
\text{and} \ H^q(Y, \Omega^p)\neq 0 \ \text{for \ some \ } p\ge \kappa, p\neq q
\notag\\
\text{if \ and \ only \ if \ } \kappa d\le n.\notag }
Let us write $y_1f_1+\ldots + y_\kappa f_\kappa$ for the equation of
$Y$, where $y_i$ are
the homogeneous coordinates in $\P^{\kappa -1}$ and $f_i$
are homogeneous polynomials of degree $d$. Then $Y$ is smooth if and only if
the codimension $\kappa$ subvariety $X\subset \P^n$ defined by
$f_1=\ldots=f_\kappa=0$ is smooth. Moreover, $Y$ is the blow-up of
$X$ in $\P^n$. Consequently, one has $CH_i(Y)=CH_{i-\kappa+1}(X)
\oplus
CH^{n+\kappa -2-i}(\P^n)$. Thus
$CH_i(Y)\otimes \Q=\Q$ for $i\le (\kappa -2)$ and
$CH_0(X)=\Z$ implies
$CH_{\kappa -1}(Y)\otimes \Q=
\Q\oplus \Q$.

Roitman's theorem is a special case of a deeper theorem for
abstractly defined Fano varieties, due to Campana, and
Koll\'ar-Miyaoka-Mori (\cite{Ko}). It is  stemming from geometry.
If $X$ is a Fano variety, then it is rationally connected, that is
 any two closed points are linked by a chain of rational curves.
This implies $CH_0(X)=\Z$, but is stronger than this. For example, in
characteristic 0, surfaces with $H^m(X, \sO_X)=0, m=1,2,$ but with nonnegative
Kodaira dimension $\le 1$ have $CH_0(X)=\Z$ (\cite{BKL}),
and certainly they are not rationally connected.
As recalled in the abstract, Fano means that $\omega_X^\vee$ is ample. Thus this strong negativity condition on the top differential forms
implies rational connectivity. On the other hand,
 strong negativity on
the 1-forms implies rationality: by the fundamental theorem of
 Mori (\cite{Mo}), $(\Omega^1_X)^\vee$ ample is equivalent to $X$ being
isomorphic to the projective space. Thus in this case $CH_0(X)=\ldots =
CH_d(X)=\Z$ and $\kappa=(d+1)$.

In consequence, it is tempting to think that the condition
\eqref{11} might result from a strong negativity condition on $d,
(d-1),\ldots, (d-\kappa+1)$ differential forms. We have at
disposal Demailly's positivity notion \cite{Dem}. A vector bundle
$E$ on a smooth projective complex variety is $s$-positive if its
hermitian curvature form, seen on $T_X\otimes E$, is positive on
all tensors of length $\le s$. Demailly's vanishing theorem says
that if $E$ is $s$ positive, then $H^q(X, \omega_X\otimes E)=0$
for $q\ge d-s+1$. In particular, let us assume that \ga{13}{
(\Omega^{n-p})^\vee \ \text{is \ } (n-p)-\rm{ positive \ for \ }
0\le p \le (\kappa -1). } Then $H^q(X, \Omega^p_X)=0 $ for $q>p,
0\le p\le (\kappa -1)$. But by Hodge duality, this implies
$H^q(X, \Omega^p_X)=0$ for $q<p, 0\le p\le (\kappa -1)$ as well.
Thus under the assumption \eqref{13}, the Beilinson-Bloch
conjectures predict \eqref{11}. One may hope that this geometric
formulation yields more information, as we discussed  for $p=0$
and $p= (d-1)$ above. It is further to be remarked that, while
applied to smooth complete intersections, Demailly's positivity
is stronger than what  would be needed to prove exactly
\eqref{11}, which is coming from the Hodge type. It is then
likely that a  positivity notion, a bit weaker than Demailly's
one,
 will force Demailly's vanishing, and would
be such that while applied to smooth complete intersections, it
would yield, via the Beilinson-Bloch conjectures, exactly the
right predicted
 statement
on Chow groups.

\section{Singular projective varieties}

In this section, $X\subset \P^n$ is still projective, but no longer
necessarily smooth.

Let us assume first that $X$ is defined over a finite field $\F_q$.
We still have that the divisibility of the eigenvalues of Frobenius
acting on $H^i_c(\bar{U},\Q_\ell)$ implies the divisibility of
 the
reciprocal roots and poles of the $\zeta$ function via
\eqref{6} and \eqref{7}. But the converse is not a priori clear.
The problem is that
the absolute values of the eigenvalues of $F_i$ are no longer determined
by $i$, there might be some
cancellation in \eqref{7}.

Next
we think of
$H^i_c(\bar{U}, \Q_\ell)$ as being no longer a semisimple
Galois representation, and
we consider its associated graded semisimple Galois representation.
 We apply the
mechanism explained in the smooth case to predict that the Hodge type of
$X$ in characteric 0 should be the same as the $\kappa$ such that $q^\kappa$
divides the eigenvalues of the Frobenius acting on $H^i_c(\bar{U}, \Q_\ell)$.
Now we have only one class of examples: projective varieties
$X$ defined by $r$ equations of degrees $d_1\ge d_2\ge \ldots \ge d_r$.
We don't require smoothness, not even that the dimension of $X$ be $(n-r)$. If
$\frac{n-d_2-\ldots -d_r}{d_1}\ge 0$ we define $\kappa$ as in \eqref{10}.
Then the theorem of Ax and Katz (\cite{Ka})
asserts that \eqref{3} is true, while
\cite{DD}, \cite{E}, \cite{ENS} show that the Hodge type of $X$ is $\kappa$.
Those bounds are sharp. In particular, we see that the coincidence of the
$\kappa$ coming from the $\zeta$
function and the one coming from
Hodge theory in this mixed case predicts
that the eigenvalues of Frobenius will be
divisible by $q^\kappa$, which we don't know so far.

In the same range of ideas, if we require now that $X$ be a
(nonsmooth) complete intersection, by \cite{W} one has higher
divisibility of the reciprocal poles of $\zeta$ (suitably
normalized), while by \cite{EW} one has a better Hodge type for
all the primitive cohomology beyond the middle dimension. Here
one would also expect that the better divisibility is not only for
the reciprocal poles of the $\zeta$ function, but also for the
eigenvalues of Frobenius acting on the corresponding cohomology
$H^i_c(\bar{U}, \Q_\ell)$ for all $i$ modulo 2 corresponding to
the poles. But in addition, in light of the Hodge type
computation, one would expect the better divisibility beyond the
middle dimension.

The Beilinson-Bloch conjectures are not formulated for projective
nonsmooth varieties. It is tempting to think that motivic
cohomology in some good sense will control both the Hodge type
and the congruence (in particular the existence) of rational
points over a finite field. However, in absence of a clear view
of what would correspond to the easy implication in the smooth
case (trivial Chow groups implies congruence for points and
nontrivial Hodge type), it is hard to forsee a good formulation
of what would be the Beilinson-Bloch conjectures in the projective
singular case.
\\ \ \\

\noindent {\it Acknowledgements}. We heartily thank
 Jean-Pierre Demailly for
showing us his concept of positivity. We  thank the participants
of the seminar at the \'Ecole Normale Sup\'erieure, in particular
Jean-Louis Colliot-Th\'el\`ene and
Ofer Gabber, for their helpful questions and comments. We also thank
Marc Levine, V. Srinivas and Daqing Wan for interesting  discussions.

\bibliographystyle{plain}

\renewcommand\refname{References}

\end{document}